\documentclass{amsart}

\usepackage{macros}

\begin{document}

\title[On the distribution of the major index on standard Young tableaux]{On the distribution of the major index \\ on standard Young tableaux}

\author{Sara C. Billey, Matja\v z Konvalinka, \and Joshua P. Swanson}

\thanks{The first author was partially supported by the Washington Research Foundation and  DMS-1764012. The second author was partially supported by Research Project BI-US/16-17-042 of the Slovenian Research Agency and research core funding No. P1-0294.}

\address{
Department of Mathematics, University of Washington, Seattle, WA 98195, USA \\
Faculty of Mathematics and Physics, University of Ljubljana \& Institute of Mathematics, Physics and Mechanics, Ljubljana, Slovenia \\
Department of Mathematics, University of California, San Diego (UCSD), La Jolla, CA  92093-0112, USA}



\begin{abstract}
The study of permutation and partition statistics is a classical topic
in enumerative combinatorics. The major index
statistic on permutations was introduced a century ago by Percy MacMahon in his
seminal works. In this extended abstract, we study the well-known generalization of the major index
to standard Young tableaux. We present several new results. In one
direction, we introduce and study two partial orders on the standard Young tableaux
  of a given partition shape, in analogy with the strong and weak
  Bruhat orders on permutations.  The existence of such
  ranked poset structures allows us to classify the realizable major
  index statistics on standard tableaux of arbitrary straight shape
  and certain skew shapes, and has representation-theoretic consequences,
  both for the symmetric group and for Shephard--Todd groups. In a different
  direction, we consider the distribution of the major index on
  standard tableaux of arbitrary straight shape and certain skew
  shapes. We classify all possible limit laws
  for any sequence of such shapes in terms of a simple auxiliary
  statistic, $\aft$, generalizing earlier results of
  Canfield--Janson--Zeilberger, Chen--Wang--Wang, and others.
  We also study unimodality, log-concavity, and local limit properties.
\end{abstract}





\maketitle

\section{Introduction}

For a skew partition $\lambda/\mu \vdash n$, denote by $\SYT(\lambda/\mu)$
the set of all standard Young tableaux of skew shape $\lambda/\mu$, i.e.~the set of all fillings of the cells of the diagram of $\lambda/\mu$ with integers $1,\ldots,n$ that
are increasing in rows and columns. We say $i$ is a \emph{descent} in a
standard tableau $T$ if $i+1$ appears in a lower row
in $T$ than $i$, where we draw partitions in English notation.
Let $\maj(T)$ denote the \emph{major index statistic} on
$\SYT(\lambda/\mu)$, which is defined to be the sum of the descents of $T$.

This statistic is a generalization of the major index on permutations or words, defined by
MacMahon in the early 1900's \cite{MacMahon.1913}
as the sum of all $i$ for which $\pi_{i+1} > \pi_i$.
The distribution of the major
index on words is a classic and surprisingly deep topic, implicitly going back before
MacMahon to Sylvester's 1878 proof of the unimodality of the $q$-binomial coefficients
\cite{doi:10.1080/14786447808639408} and
beyond\footnote{Indeed, Sylvester's excitement at settling
this then-quarter-century-old conjecture is palpable: ``I accomplished with scarcely an effort
a task which I had believed lay outside the range of human power.''}.
In this extended abstract, we summarize several recent explorations of the distribution of the
major index on ordinary shape tableaux and certain skew tableaux.

The major index generating function for $\SYT(\lambda/\mu)$ for a straight shape
$\lambda/\mu$ is given by
\begin{equation}\label{eq:syt-q}
  \SYT(\lambda/\mu)^{\maj}(q) := \sum_{T \in \SYT(\lambda/\mu)} q^{\maj(T)}
    = \sum_{k\geq 0} b_{\lambda/\mu,k} q^k
\end{equation}
where the coefficients are the \emph{fake degree sequence}
\begin{equation}\label{eq:a_la}
  b_{\lambda/\mu, k} := \#\{T \in \SYT(\lambda/\mu) : \maj(T) = k\} \text{ for }  k=0,1,2,\ldots
\end{equation}
The fake degrees for straight shapes $\lambda$ and certain skew shapes
have appeared in a variety of algebraic and
representation-theoretic contexts  including Green's work on the
degree polynomials of unipotent $\GL_n(\bF_q)$-rep\-re\-sen\-ta\-tions
\cite[Lem.~7.4]{MR0072878}, the irreducible decomposition of type $A$
coinvariant algebras \cite[Prop.~4.11]{stanley.1979}, Lusztig's work
on the irreducible representations of classical groups
\cite{Lusztig.77}, and branching rules between symmetric groups and
cyclic subgroups \cite[Thm.~3.3]{stembridge89}.  The term ``fake
degree'' was apparently coined by Lusztig \cite{Carter.89}, most likely
because $\#\SYT(\lambda)= \sum_{k\geq 0} b_{\lambda,k}$ is the degree
of the irreducible $S_n$-representation indexed by $\lambda$, so a
$q$-analog of this number is not itself a degree but is related to the degree.

We consider three natural enumerative questions involving the fake degrees:
\begin{enumerate}
  \item[Q1.] Which $b_{\lambda,k}$ are zero?
  \item[Q2.] Are there efficient asymptotic estimates for $b_{\lambda, k}$?
  \item[Q3.] Are the fake degree sequences unimodal?
\end{enumerate}

Our results are presented in full in \cite{bks2, bks1}. Given the length of these papers and space limitations for this extended abstract, we present no proofs, and instead provide very rough sketches with references to the full papers. We describe the answer to Q1 in
\Cref{sec:zeros}, a complete answer to one precise version of Q2 in \Cref{sec:an},
and further work and open problems related to Q3 and beyond in \Cref{sec:more}.

\section{Zeros of the fake degree sequence}\label{sec:zeros}

We completely settle Q1 with the following result. Write $b(\lambda)
:= \sum (i-1)\lambda_i$ and let $\lambda'$ denote the conjugate partition of $\lambda$.

\begin{Theorem}\label{thm:zeros}
  For every partition $\lambda\vdash n\geq1$ and integer $k$ such that
  $b(\lambda) \leq k \leq \binom{n}{2}-b(\lambda')$, we have
  $b_{\lambda, k} > 0$ except in the case when $\lambda$ is
  a rectangle with at least two rows and columns and $k$ is either
  $b(\lambda)+1$ or $\binom{n}{2}-b(\lambda')-1$. Furthermore,
  $b_{\lambda, k} = 0$ for $k < b(\lambda)$ or $k > \binom{n}{2}
  - b(\lambda')$.
\end{Theorem}

The main ingredient of the proof is a map $\varphi \, \colon \SYT(\lambda) \setminus \cE(\lambda) \rightarrow \SYT(\lambda)$
with the property $\maj(\varphi(T)) = \maj(T) + 1$. Here $\cE(\lambda)$ is the (small) set of exceptional tableaux where such a map cannot be defined and contains:
\begin{enumerate}
    \item[i.] For all $\lambda$, the tableau for $\lambda$ with the largest possible major index.

    \item[ii.] If $\lambda$ is a rectangle, the tableau for $\lambda$ with the smallest possible major index.

    \item[iii.] If $\lambda$ is a rectangle with at least two rows
      and columns, the unique tableau in $\SYT(\lambda)$ with
      major index equal to $\binom{n}{2}-b(\lambda')-2$.
  \end{enumerate}

The construction of such a map goes roughly as follows. For most tableaux, we can apply a simple \emph{rotation rule} that increases the major index by $1$. More specifically, given $T \in \SYT(\lambda)$, assume that we have an interval $[i, k] \subset [n]$ such that $T' := (i, i+1, \ldots, k-1, k) \cdot T$ is in $\SYT(\lambda)$ and such that there is some $j$ for which
    \[ \{j\} = \Des(T') - \Des(T)\qquad\text{ and }\qquad
        \{j-1\} = \Des(T) - \Des(T'). \]
For example, assume $i < j < k$. Then we require that:
\begin{itemize}
  \item $i, \ldots, j-1$ form a horizontal strip, $j-1, j$ form a vertical
      strip, and $j, j+1, \ldots, k$ form a horizontal strip;
  \item $i$ appears strictly northeast of $k$ and
      $i-1$ is not in the rectangle bounding $i$ and $k$;
  \item $k$ appears strictly northeast of $k-1$ and
      $k+1$ is not in the rectangle bounding $k$ and $k-1$;
\end{itemize}

\begin{figure}[ht]
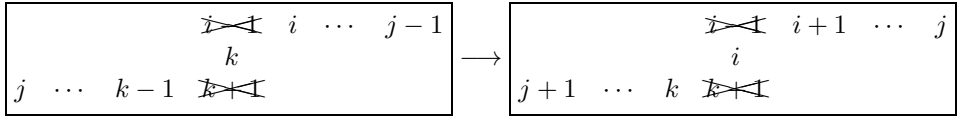

  \centering
    \[
      \boxed{\begin{matrix}
        & & & \xcancel{i-1} & i &\cdots & j-1 \\
        & & & k \\
        j &\cdots & k-1 & \xcancel{k+1} \\
      \end{matrix}}
      \longrightarrow
      \boxed{\begin{matrix}
        & & & \xcancel{i-1} & i+1 & \cdots & j \\
        & & & i \\
        j+1 & \cdots & k & \xcancel{k+1}\\
      \end{matrix}}
    \]
    \caption{A positive rotation with $i < j < k$.}
  \label{fig:posrot}
\end{figure}

A sketch of the positive rotation in this case is presented in Figure \ref{fig:posrot}. The cases $i = j$ and $j = k$ are slightly different. If $i, j, k$ with the required properties exist, we take the lexicographically smallest such numbers, and define $\varphi(T) := T'$.

We similarly define \emph{negative rotations} where $T' = (k, k-1, \ldots, i) \cdot T$. If no positive rotations are possible on a tableau $T$, but we have negative rotations, use the lexi\-co\-graphically smallest one to define $\varphi(T)$.

There are very few cases when neither a positive nor a negative rotation can be applied. For example, among the $81,081$ tableaux in
$\SYT(5442)$, there are only $24$ on which we cannot
apply any positive or negative rotation rule.  In particular, no
rotation rules can be applied to the following two tableaux:
\[
  \begin{matrix}
    1 & 2 & 3 & 4 & 5 \\
    6 & 7 & 8 & 9 & \text{} \\
    10 & 11 & 12 & 13 & \text{} \\
    14 & 15 & \text{} & \text{} & \text{} \\
  \end{matrix} \qquad \text{and} \qquad
  \begin{matrix}
    1 & 2 & 3 & 8 & 12 \\
    4 & 6 & 9 & 13 & \text{} \\
    5 & 7 & 10 & 14 & \text{} \\
    11 & 15 & \text{} & \text{} & \text{} \\
  \end{matrix}.
\]
In that case, we have to apply what we call \emph{block rules} B1--B5. We refer the reader to \cite[Def.~4.13]{bks1} for explicit definitions. 



As a consequence of the proof of \Cref{thm:zeros}, we identify two new
ranked poset structures on $\SYT(\lambda)$ where the rank function is
determined by $\maj$.  Furthermore, as a corollary of \Cref{thm:zeros}
we have a new proof of a complete classification due to the third
author \cite[Thm.~1.4]{s17} generalizing an earlier result of Klyachko
\cite{klyachko74} for when the counts
\[ a_{\lambda, r} := \{T \in \SYT(\lambda) : \maj(T) \equiv_n r\}
\text{ for } \lambda \vdash n \]
are nonzero.
We can also classify internal zeros of certain skew shapes,
\cite[Lem.~6.2]{bks1}.

\begin{Remark}
Lascoux--Sch\"utzenberger \cite{LS7} defined an operation called
\textit{cyclage} on semi\-standard tableaux, which decreases
\textit{cocharge} by $1$. The \textit{cyclage poset} on the set of
semi\-standard tableaux arises from applying cyclage in all possible
ways. Cyclage preserves the \textit{content}, i.e.~the number of
$1$'s, $2$'s, etc.  See also \cite[Sect. 4.2]{Shimozono-Weyman.2000}.
Restricting to standard tableaux, cocharge coincides with $\maj$, so
the cyclage poset on $\SYT(n)$ is ranked by $\maj$. However, cyclage
does not necessarily preserve the shape, so it does not suffice to
prove \Cref{thm:zeros}.  For example, restricting the cyclage poset to
$\SYT(32)$ gives a poset which has two connected components and is
not ranked by $\maj$, while both of our poset structures on
$\SYT(32)$ are chains.  One reviewer posed an interesting question: is
there any relation between the cyclage poset covering relations
restricted to $\SYT(\lambda)$ and the two ranked poset structures used
to prove \Cref{thm:zeros}?  We have not found one.
\end{Remark}

Symmetric groups are the finite reflection groups of type $A$.  The
classification and invariant theory of both finite irreducible real
reflection groups and complex reflection groups developed over the
past century builds on our understanding of the type $A$ case
\cite{Hum}.  In particular, these groups are classified by
Shephard--Todd into an infinite family $G(m,d,n)$ together with $34$ exceptions.
Using work of Stembridge on generalized exponents for irreducible
representations, the analog of \eqref{eq:syt-q}
can be phrased for all Shephard--Todd groups as
\begin{equation}\label{eq:gmdn}
    \gmdn{\lambda}(q) :=
      \frac{\#\{\underline{\lambda}\}^d}{d} \cdot
      \bbinom{n}{\alpha(\underline{\lambda})}_{q; d} \cdot
      \prod_{i=1}^m \SYT(\lambda^{(i)})^{\maj}(q^m) =
      \sum b_{\{\underline{\lambda}\}^d, k} q^k
    \end{equation}
    where
    $\underline{\lambda}=(\lambda^{(1)}, \ldots, \lambda^{(m)})$ is
    a sequence of $m$ partitions with $n$ cells total,
    $\alpha(\underline{\lambda})=(|\lambda^{(1)}|,
    \ldots,|\lambda^{(m)}|)\vDash n$, $d \mid m$, and
    $\{\underline{\lambda}\}^d$ is the orbit of $\underline{\lambda}$
    under the group $C_d$ of $(m/d)$-fold cyclic rotations;
    see \cite[Cor.~8.2]{bks1}.
The polynomials
    $ \bbinom{n}{\alpha(\underline{\lambda})}_{q; d} $ are
    deformations of the usual $q$-multinomial coefficients which we
    explore in \cite[\S 7]{bks1}.
    The
    coefficients $b_{\{\underline{\lambda}\}^d, k}$ are the fake
    degrees in this case.

    We use \eqref{eq:gmdn} and \Cref{thm:zeros} to completely classify
    all nonzero fake degrees for coinvariant algebras for all
    Shephard--Todd groups $G(m,d,n)$, which includes the finite
    real reflection groups in types $A$, $B$, and $D$.  See \cite[\S 8]{bks1} for details.

\section{Asymptotic normality of the major index on SYT}\label{sec:an}

Let us turn our attention to question Q2. It is a well-known fact that the major index
statistic on permutations satisfies a central limit theorem.
Given a real-valued random variable $\cX$, we let $\cX^* := \frac{\cX - \mu}{\sigma}$
denote the corresponding normalized random variable with mean $0$
and variance $1$.

\begin{Theorem}
  \cite{MR0013252}
  \label{thm:maj_Sn_an}
  Let $\cX_n[\maj]$ denote the major index random variable on $S_n$ under the
  uniform distribution. Then, for all $t \in \bR$,
    \[ \lim_{n \to \infty} \bP[\cX_n[\maj]^* \leq t] = \bP[\cN \leq t] \]
  where $\cN$ is the standard normal random variable.
\end{Theorem}

\begin{figure}[ht]
  \centering
  \begin{subfigure}[t]{0.3\textwidth}
    \centering
    \includegraphics[width=\textwidth]{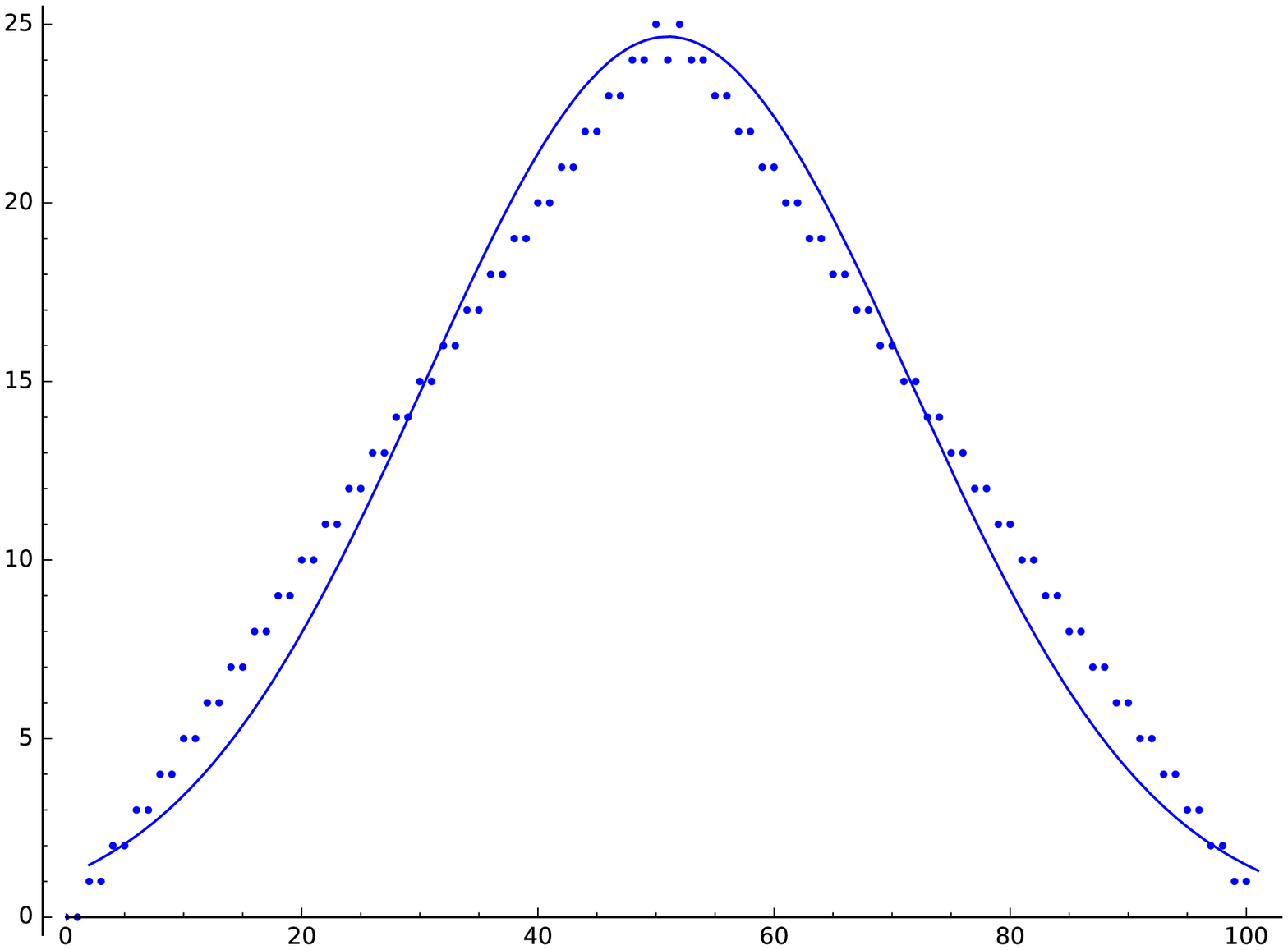}
    \caption{$\lambda = (50, 2)$, $\aft(\lambda) = 2$}
    \label{fig:an_examples_a}
  \end{subfigure}
  \begin{subfigure}[t]{0.3\textwidth}
    \centering
    \includegraphics[width=\textwidth]{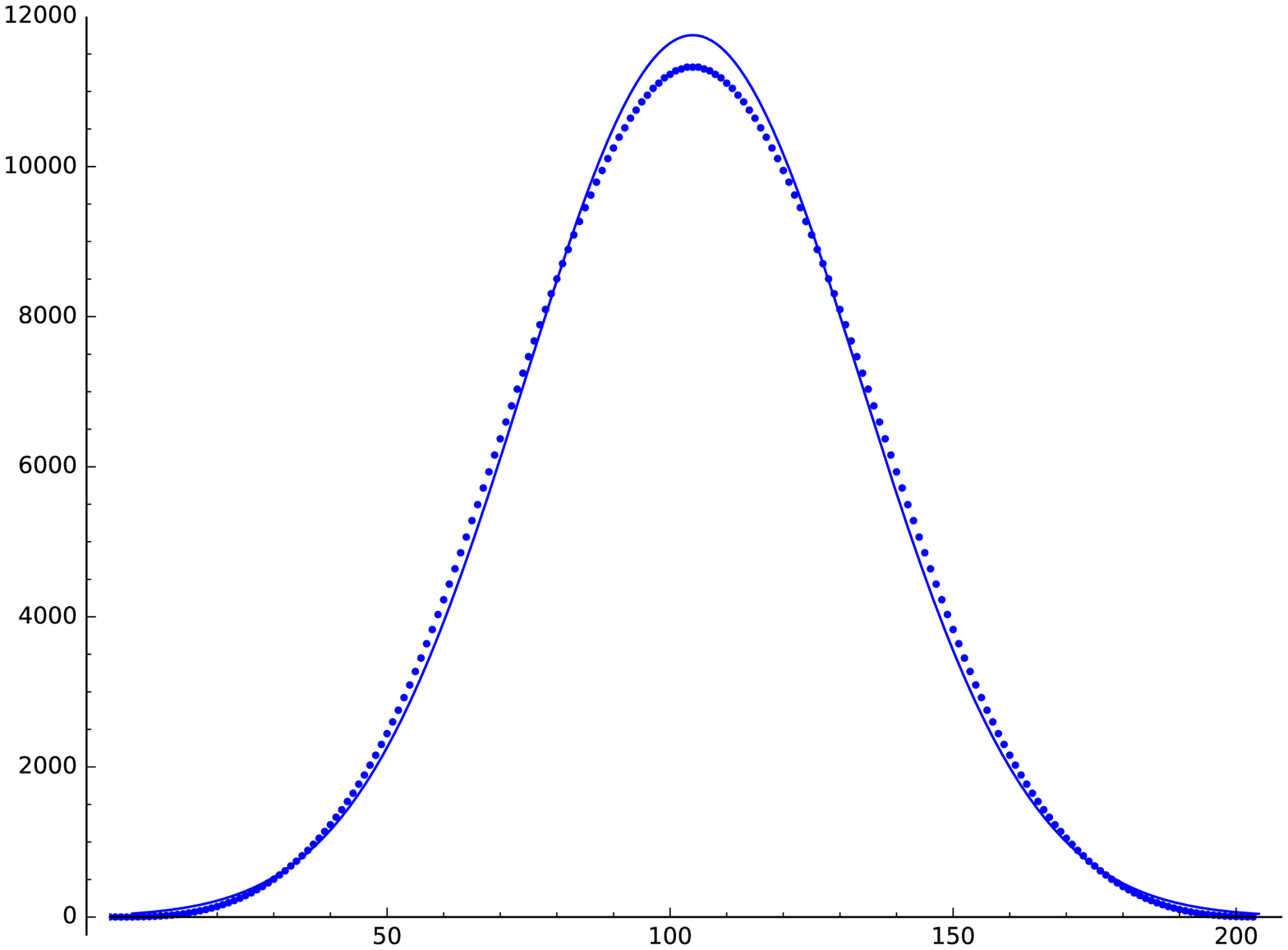}
    \caption{$\lambda = (50, 3, 1)$, $\aft(\lambda) = 4$}
    \label{fig:an_examples_b}
  \end{subfigure}
  \begin{subfigure}[t]{0.35\textwidth}
    \centering
    \includegraphics[width=\textwidth]{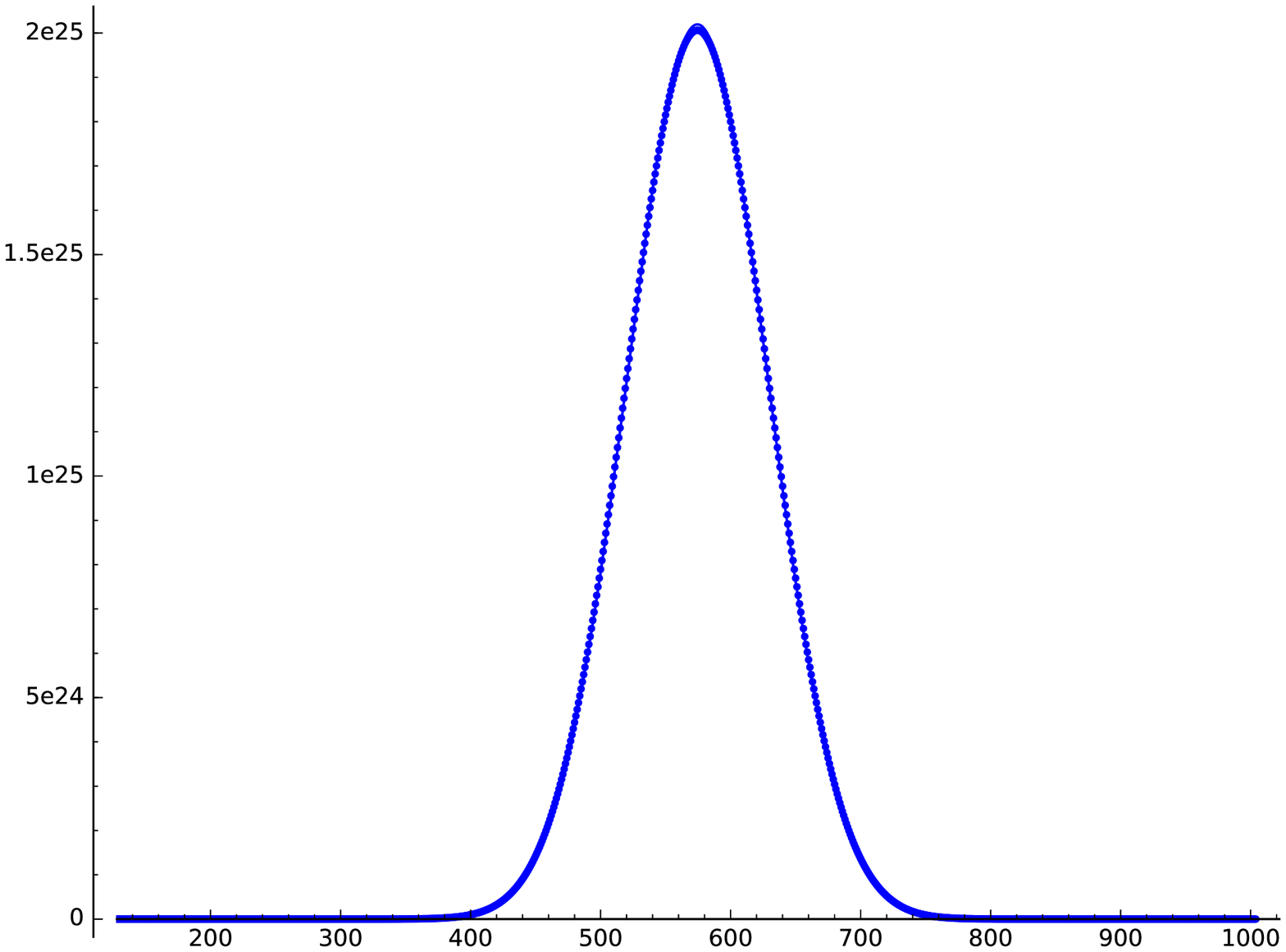}
    \caption{$\lambda = (8^{2}, 7, 6, 5^{3}, 2^{2})$, $\aft = 39$}
    \label{fig:an_examples_c}
  \end{subfigure}
  \caption{Plots of $\#\{T \in \SYT(\lambda) : \maj(T) = k\}$
     as a function of $k$ for three partitions $\lambda$, overlaid
     with scaled Gaussian approximations using the same
     mean and variance.}
  \label{fig:an_examples}
\end{figure}

Briefly, we say $\maj$ on $S_n$ is
\textit{asymptotically normal} as $n \to \infty$. See \cite[Table 1]{bks2} for further examples of asymptotic normality.
\Cref{fig:an_examples} shows some sample distributions for the
major index on standard tableaux for three particular partition shapes.
Note that Gaussian approximations fit the data well.

In \Cref{thm:maj_Sn_an}, we simply let $n \to \infty$. For partitions,
the shape $\lambda$ may ``go to infinity'' in many different ways.
The following statistic on partitions overcomes this difficulty.

\begin{Definition}
  Suppose $\lambda$ is a partition. Let the \textit{aft} of $\lambda$ be
    $\aft(\lambda) :=	 |\lambda| - \max\{\lambda_1, \lambda_1'\}.$
\end{Definition}
\noindent Intuitively, if the first row of $\lambda$ is at least as long as the
first column, then $\aft(\lambda)$ is the number of cells \textit{not}
in the first row. This definition is strongly reminiscent of a
\textit{representation stability} result of Church and Farb
\cite[Thm.~7.1]{MR3084430}, which is proved with an analysis of
the major index on standard tableaux.

Our first main result gives the analogue of \Cref{thm:maj_Sn_an} for
$\maj$ on $\SYT(\lambda)$. In particular, it completely classifies which
sequences of partition shapes give rise to as\-ymp\-totically normal
sequences of $\maj$ statistics on standard tableaux.

\begin{Theorem}\label{thm:an}
  Suppose $\lambda^{(1)}, \lambda^{(2)}, \ldots$ is a sequence of
  partitions, and let $\cX_N = \cX_{\lambda^{(N)}}[\maj]$ be the corresponding
  random variables for the $\maj$ statistic on $\SYT(\lambda^{(N)})$.
  Then, the sequence $\cX_1, \cX_2, \ldots$ is asymptotically normal if and only if
  $\aft(\lambda^{(N)}) \to \infty$ as $N \to \infty$.
\end{Theorem}

\begin{Remark}
  In \cite{bks2}, we more generally consider $\maj$
  on $\SYT(\underline{\lambda})$ where $\underline{\lambda}$ is a
  \textit{block diagonal} skew partition. Special cases of this include
  Canfield--Janson--Zeilberger's main result in \cite{MR2794017}
  classifying asymptotic normality for $\inv$ or $\maj$ on words
  (though see \cite{MR2925927c} for earlier, essentially equivalent
  results due to Diaconis \cite{MR964069}).
  The case of words generalizes
  \Cref{thm:maj_Sn_an}. The $\lambda^{(N)} = (N, N)$ case of
  \Cref{thm:an} also recovers the main result of Chen--Wang--Wang
  \cite{Chen-Wang-Wang.2008}, giving asymptotic normality for
  $q$-Catalan coefficients.
\end{Remark}

Our proof of \Cref{thm:an} relies on the \textit{method of moments},
which requires useful descriptions of the moments of $\cX_\lambda[\maj]$.
Adin--Roichman \cite{MR1841639} gave exact formulas for the mean
and variance of $\cX_\lambda[\maj]$ in terms of the
{hook lengths} of $\lambda$. These formulas are obtained from Stanley's
elegant closed form for the polynomials $\SYT(\lambda)^{\maj}(q)$.
Let $h_c = \lambda_i + \lambda'_j - i - j + 1$ denote the hook
length of the cell $c = (i,j)$.

\begin{Theorem}{\cite[7.21.5]{ec2}}\label{thm:stanley_maj}
  Let $\lambda \vdash n$ with $\lambda = (\lambda_1, \lambda_2, \ldots)$. Then
  \begin{equation}\label{eq:stanley_maj}
    \SYT(\lambda)^{\maj}(q)
       = \frac{q^{b(\lambda)}[n]_q!}
                  {\prod_{c \in \lambda} [h_c]_q}.
  \end{equation}
\end{Theorem}

More generally, formulas for the $d$th moment $\mu_d^\lambda$,
$d$th central moment $\alpha_d^\lambda$, and
$d$th \textit{cumulant} $\kappa_d^\lambda$ of $\maj$ on
$\SYT(\lambda)$ may be derived from \Cref{thm:stanley_maj}. Here the cumulants
$\kappa_1, \kappa_2, \ldots$ of $\cX$ are defined
  to be the coefficients of the exponential generating function
    \[ K_{\cX}(t) := \sum_{d=1}^\infty \kappa_d \frac{t^d}{d!}
       := \log M_{\cX}(t) = \log \bE[e^{t\cX}]. \]
The most elegant of these
formulas is for the cumulants, from which the moments and
central moments are all easy to compute.
\begin{Theorem}\label{thm:SYT_moments}
  Let $\lambda \vdash n$ and $d \in \bZ_{> 1}$. We have
  \begin{equation}\label{eq:SYT_cumulants}
    \kappa_d^\lambda =  \frac{B_d}{d}
    \left[ \sum_{j=1}^n j^d - \sum_{c \in \lambda} h_c^d \right]
  \end{equation}
  where $B_0, B_1, B_2, \ldots = 1, \frac{1}{2}, \frac{1}{6}, 0, -\frac{1}{30}, 0, \frac{1}{42},
  0, \ldots$ are the Bernoulli numbers.
\end{Theorem}
\noindent See \cite[Thm.~2.9]{bks2} for a
generalization of \eqref{eq:SYT_cumulants}
along with exact formulas for the moments and
central moments. See \cite[Rem.~2.10]{bks2} for
some of the history of this formula.

  For ``most'' partition shapes, one expects the term $\sum_{j=1}^n j^d$
  in \eqref{eq:SYT_cumulants} to dominate $\sum_{c \in \lambda} h_c^d$,
  in which case asymptotic normality is quite straightforward. However,
  for some shapes there is a very large amount of cancellation in
  \eqref{eq:SYT_cumulants} and determining the limit law can be quite
  subtle.

  \begin{Remark}
  $\cX_{\lambda}[\maj]$ can be written as
  the sum of scaled indicator random variables $D_1, 2D_2,3D_3, \ldots,$ $(n-1)D_{n-1}$
  where $D_i$ determines if there is a descent at position $i$. However, the $D_i$ are not at all independent,
  so one may not simply apply standard central limit theorems.
  Interestingly, the $D_i$ are identically distributed \cite[Prop.~7.19.9]{ec2}.
  The lack of independence of the $D_i$'s likewise complicates related work
  by Fulman \cite{MR1652841} and Kim--Lee \cite{KIM2020105123}
  considering the limiting distribution of descents.
\end{Remark}

The non-normal continuous limit laws for $\maj$ on $\SYT(\lambda)$ turn out
to be the \textit{Irwin--Hall distributions} $\cIH_M := \sum_{k=1}^M
\cU[0, 1]$, which are the sum of $M$ i.i.d.~continuous $[0, 1]$ random variables.
The following result completely classifies all possible limit laws for
$\maj$ on $\SYT(\lambda)$ for any sequence of partition shapes.

\begin{Theorem}\label{thm:all_limits}
  Let $\lambda^{(1)}, \lambda^{(2)}, \ldots$ be a sequence of partitions.
  Then $(\cX_{\lambda^{(N)}}^*[\maj])$ converges in distribution if and only if
  \begin{enumerate}
    \item $\aft(\lambda^{(N)}) \to \infty$; or
    \item $|\lambda^{(N)}| \to \infty$ and $\aft(\lambda^{(N)}) \to M < \infty$; or
    \item the distribution of $\cX_{\lambda^{(N)}}^*[\maj]$ is eventually constant.
  \end{enumerate}
  The limit law is $\cN$ in case (i), $\cIH_M^*$ in
  case (ii), and discrete in case (iii).
\end{Theorem}

Case (iii) naturally leads to the question: when does $\cX_{\lambda}^*[\maj]
= \cX_{\mu}^*[\maj]$? Such a description in terms of hook lengths is
given in \cite[Thm.~7.1]{bks2}.

\begin{Example}
  We illustrate each possible limit. For (i), let
  $\lambda^{(N)} := (N, \lfloor\ln N\rfloor)$, so that
  $\aft(\lambda^{(N)}) = \lfloor\ln N\rfloor \to \infty$ and the
  distributions are asymptotically normal.
  For (ii), fix $M \in \bZ_{\geq 0}$ and let $\lambda^{(N)}
  := (N+M, M)$, so that $\aft(\lambda^{(N)}) = M$ is constant
  and the distributions converge to $\cIH_M^*$.
  For (iii), let $\lambda^{(2N)} := (12, 12, 3, 3, 3, 2, 2, 1, 1)$
  and $\lambda^{(2N+1)} := (15, 6, 6, 6, 4, 2)$, which have the
  same multisets of hook lengths despite not being transposes of each other, and
  consequently the same normalized
  distributions.
\end{Example}

In order to be able to use the method of moments, we need the following  result.

\begin{Theorem}[Frech\'et--Shohat Theorem,
  {\cite[Thm.~30.2]{MR1324786}}]\label{thm:moments}
  Let $\cX_1, \cX_2, \ldots$ be a sequence of real-valued random
  variables,
  and let $\cX$ be a real-valued random variable. Suppose
  the moments of $\cX_n$ and $\cX$ all exist and the moment
  generating functions all have a positive radius of convergence. If
  \begin{equation}\label{eq:moments_criterion}
    \lim_{n \to \infty} \mu_d^{\cX_n} = \mu_d^{\cX} \hspace{.5cm} \forall
     d \in \bZ_{\geq 1},
  \end{equation}
  then $\cX_1, \cX_2, \ldots$ converges in distribution to $\cX$.
\end{Theorem}

By \Cref{thm:moments} we may test for asymptotic normality on
level of individual nor\-mal\-ized moments, which is often referred to as the
\textit{method of moments}. By the formula
\begin{align}\label{eq:moment_to_cumulant}
  \mu_d = \kappa_d + \sum_{m=1}^{d-1}
    \binom{d-1}{m-1} \kappa_m \mu_{d-m},
\end{align}
which is not hard to derive, we may further replace the moment condition
\eqref{eq:moments_criterion} with the corresponding cumulant
condition. For instance, we have the following explicit criterion.

\begin{Corollary}\label{cor:cumulants}
  A sequence $\cX_1, \cX_2, \ldots$ of real-valued
  random variables on finite sets is as\-ymp\-totically normal
  if for all $d \geq 3$ we have
  \begin{equation}\label{eq:cumulants_criterion2}
    \lim_{n \to \infty} \frac{\kappa^{\cX_n}_d}{(\sigma^{\cX_n})^d} = 0.
  \end{equation}
\end{Corollary}

The hardest part of the proof of \Cref{thm:all_limits} is asymptotic normality. We can prove that part by using the following three lemmas. We refer to \cite[\S 5]{bks2} for proofs and further background.

\begin{Definition}
  A \textit{reverse standard Young tableau} of shape $\lambda/\nu$
  is a bijective filling of $\lambda/\nu$ which strictly
  decreases along rows and columns. The set of reverse
  standard Young tableaux of shape $\lambda/\nu$ is denoted
  $\RSYT(\lambda/\nu)$.
\end{Definition}

\begin{Lemma}\label{lem:rsyt}
  Let $\lambda/\nu \vdash n$ and $T \in \RSYT(\lambda/\nu)$. Then
  for all $c \in \lambda/\nu$,
  \begin{align}\label{eq:rsyt1}
    T_c \geq h_c.
  \end{align}
  Furthermore, for any positive integer $d$,
  \begin{align}\label{eq:rsyt2}
    \sum_{j=1}^n j^d - \sum_{c \in \lambda/\nu} h_c^d
      &= \sum_{c \in \lambda/\nu} (T_c^d - h_c^d)
        = \sum_{c \in \lambda/\nu} (T_c - h_c) {\mathbf{h}}_{d-1}(T_c, h_c),
  \end{align}
  where ${\mathbf h}_{d-1}$ denotes the complete homogeneous symmetric function.
\end{Lemma}

\begin{Lemma}\label{lem:small_hook_bounds}
  Let $\lambda/\nu \vdash n$ such that
  $\max_{c \in \lambda/\nu} h_c < 0.8n$. Let $d$ be
  any positive integer. Then
  \[ \frac{n^{d+1}}{26(d+1)} - 2(0.8)^d n^d
      < \sum_{j=1}^n j^d - \sum_{c \in \lambda/\nu} h_c^d
      < \frac{n^{d+1}}{d+1} + n^d. \]
\end{Lemma}

\begin{Lemma}\label{lem:large_hook_bounds}
  Let $\lambda/\nu \vdash n$ such that $\max_{c \in \lambda/\nu} h_c
  \geq 0.8n$, and let $d$ be any positive integer. Furthermore, suppose
  $n \geq 10$. Then,
  \begin{align}
     \aft(\lambda/\nu) \frac{\lfloor 0.1n\rfloor^d}{d}  \hspace{.2cm} \leq  \hspace{.2cm}
    \sum_{j=1}^n j^d - \sum_{c \in \lambda/\nu} h_c^d  \hspace{.2cm}
    \leq  \hspace{.2cm} 2\aft(\lambda/\nu) \left(n^d + dn^{d-1}\right).
  \end{align}
\end{Lemma}

\begin{Corollary}\label{cor:aft_theta}
  For fixed $d \in \bZ_{\geq 1}$, uniformly for all skew shapes $\lambda/\nu$,
    \begin{equation} \sum_{k=1}^{|\lambda/\nu|} k^d - \sum_{c \in \lambda/\nu} h_c^d
      = \Theta(\aft(\lambda/\nu) \cdot |\lambda/\nu|^d). \end{equation}


  \begin{proof}
    Let $n = |\lambda/\nu|$.
    When $\max_{c \in \lambda/\nu} h_c \geq 0.8n$, the result
    follows from \Cref{lem:large_hook_bounds}. On the other hand, when
    $\max_{c \in \lambda/\nu} h_c < 0.8n$, then
    $n \geq \aft(\lambda/\nu) \geq 0.2n$, and the result follows
    from \Cref{lem:small_hook_bounds}.
  \end{proof}
\end{Corollary}

\begin{Corollary}\label{cor:bound.cumulants}
  Fix $d$ to be an even positive integer. Uniformly for
  all block diagonal skew shapes $\underline{\lambda}$,
  the absolute value of the normalized cumulant
  $|{\kappa_d^{\underline{\lambda}}}^*|$ of
  $\cX_{\underline{\lambda}}[\maj]$ is $\Theta(\aft(\underline{\lambda})^{1-d/2})$.

  \begin{proof}
    For $d$ even, by a generalization of \Cref{thm:SYT_moments} to block diagonal shapes and \Cref{cor:aft_theta},
    we have
      $ |\kappa_d^{\underline{\lambda}}|
         = \Theta(\aft(\underline{\lambda}) n^d), $
    where $n = |\underline{\lambda}|$.
    Consequently,
      \[ |{\kappa_d^{\underline{\lambda}}}^*|
          = \left|\frac{\kappa_d^{\underline{\lambda}}}
             {(\kappa_2^{\underline{\lambda}})^{d/2}}\right|
          = \Theta\left(\frac{\aft(\underline{\lambda})n^d}
             {\aft(\underline{\lambda})^{d/2} n^d}\right)
          = \Theta(\aft(\underline{\lambda})^{1-d/2})\]
 by the homogeneity of cumulants.
  \end{proof}
\end{Corollary}

A natural generalization of \Cref{thm:an} to block diagonal skew shapes
$\underline{\lambda}$ follows by combining \Cref{cor:cumulants},
\Cref{cor:bound.cumulants}, and similar estimates in the ``degenerate''
case when $\aft(\underline{\lambda})$ is bounded.  See
\cite[Thm.~6.3]{bks2}.

\section{Unimodality and beyond}\label{sec:more}

The easy answer to question Q3 is ``no'':  the fake degree sequences
are not always unimodal.  For example, $ \SYT(42)^{\maj}(q)$ is not
unimodal. Nonetheless, certain inversion
number generating functions $p_\alpha^{(k)}(q)$ which appear in
a generalization of $\SYT(\lambda)^{\maj}(q)$ are in
fact unimodal; see \cite[Def.~7.7, Cor.~7.10]{bks1}.
Furthermore, computational evidence suggests $\SYT(\lambda)^{\maj}(q)$
is typically not far from unimodal.  See also \cite[Sect. 8]{bks2} concerning unimodality,
log-concavity, and asymptotic normality for  skew shapes.

\begin{Conjecture}\label{conj:unimodality}
The polynomial $\SYT (\lambda )^{\maj }(q)$ is unimodal if
$\lambda $ has at least $4$ corners.
\end{Conjecture}

In another direction, one may ask for a precise estimate of the deviation of
$b_{\lambda, k}$ from normal with an explicit error bound. Such a result is called a
\textit{local limit theorem}.

\begin{Conjecture}\label{conj:local.limit}
Let $\lambda \vdash n$ be any partition. Uniformly for all $n$, for all
integers $k$, we have
\begin{equation*}
\left |\bP [X_{\lambda }[\maj ] = k] - f(k; \kappa _{1}^{\lambda },
\kappa _{2}^{\lambda })\right | = O\left (
\frac{1}{\sigma _{\lambda }\aft (\lambda )}\right ).
\end{equation*}
\end{Conjecture}

A sequence $a_0, a_1, a_2, \ldots$ is \textit{parity-unimodal} if
$a_0, a_2, a_4, \ldots$ and $a_1, a_3, a_5, \ldots$ are each
unimodal. Stucky \cite[Thm.~1.3]{stucky2018cyclic} recently showed that the
$q$-Catalan polynomials, namely $\SYT((n, n))^{maj}(q)$ up to a
$q$-shift, are parity-unimodal. The argument involves constructing an
$\fsl_2$-action on rational Cherednik algebras. See [\S3.1, Haiman94]
for a prototype of the argument in a highly related context. Based on
Stucky's result, our internal zeros classification, and a brute-force
check for $n \leq 50$, we conjecture the following.

\begin{Conjecture} The fake-degree polynomials $f^\lambda(q)$ are parity-unimodal for all $\lambda$.
\end{Conjecture}

\section{Acknowledgements}

We would like to thank Krzysztof Burdzy, Rodney Canfield, Persi
Diaconis, Sergey Fomin, Pavel Galashin, Svante Janson, William
McGovern, Alejandro Morales, Andrew Ohana, Greta Panova, Mihael Perman, Martin Rai\v{c},
Victor Reiner, Richard Stanley, Chris\-tian Stump, Sheila Sundaram, Vasu Tewari, Lauren Williams,
Alex Woo, and the referees for helpful insights related to this work.

\bibliographystyle{alpha}
\bibliography{refs}

\end{document}